\documentclass[12pt,leqno]{article}
\newcommand{\kmcomment}[1]{}

\newcommand{\ovfrakg}{\overline{\frakg}}
\newcommand{\ovfrakh}{\overline{\frakh}}

\newcommand{\ds}{\ensuremath{\displaystyle }}
\newcommand{\pdel}{\partial}

\newcommand{\myHom}[1]{\textrm{H}_{#1}}

\newcommand{\mR}{\ensuremath{\mathbb{R}}} 
\newcommand{\mZ}{\ensuremath{\mathbb{Z}}} 

\newcommand{\yb}[1]{y_{#1}}
\newcommand{\zb}[1]{z_{#1}}

\newcommand{\frakg}{\mathfrak{g}}
\newcommand{\frakgN}[1]{{\mathfrak{g}}_{#1}}

\newcommand{\frakh}{\mathfrak{h}}

\newcommand{\tbdl}[1]{\mathrm{T}(#1)}

\newcommand{\tbdli}[2]{\mathrm{T}_{#2}(#1)}
\newcommand{\cbdl}[1]{\mathrm{T}^{*}(#1)}

\newcommand{\cbdli}[2]{\mathrm{T}^{*}_{#2}(#1)}

\newcommand{\Sbt}[2]{[#1,#2]}                
                
\newcommand{\parity}[1]{(-1)^{#1}}

\newcommand{\Lb}[1]{\mathit{L}_{#1}}  
\newcommand{\inn}[1]{i_{#1}} 
\newcommand{\inner}[2]{\langle #1 , #2 \rangle}

\newcommand{\Pkt}[2]{\{#1,#2\}}%

\newcommand{\BktU}[3]{\{#2,#3\}^{#1}}

\newcommand{\Ekt}[2]{\{#1,#2\}^{'}}
\newcommand{\gaiseki}[1]{\Lambda^{#1} \tbdl{M}}

\newcommand{\cgaiseki}[1]{\Lambda^{#1} \cbdl{M}}

\newcommand{\wtedC}[2]{ C_{#1}^{[#2]}} 
\newcommand{\wtedCR}[2]{ \widetilde{C}_{#1}^{[#2]}}

\newcommand{\SbtPhi}[2]{\ensuremath{[#1,#2]^{\phi}}}

\newcommand{\mw}{\mywedge}
\newcommand{\we}{\wedge}

\newcommand{\SbtES}[2]{ [#1,#2]_{res}}

\newcommand{\dphi}{ d_{\phi}}

\renewcommand{\dim}{\textrm{dim}}

\newcommand{\mywedge}{\bigtriangleup} 

\newcommand{\myCS}[1]{ \text{C}_{#1}} 

\renewcommand{\[}{$$} \renewcommand{\]}{$$}

\usepackage{ntheorem} 
{
\theoremheaderfont{\bfseries}
\theorembodyfont{\rmfamily}

\newtheorem{defn}{\textbf{Definition}}

\newtheorem{prop}{Proposition}[section]
\newtheorem{exam}{Example}[section]
\newtheorem{remark}{Remark}[section]

\newtheorem{theorem}{Theorem}[section] 
\newtheorem{thm}{Theorem}[section]
\newtheorem*{thm-none}{Theorem}[section]

\newtheorem{kmCor}[theorem]{Corollary}

}

\newcommand{\qed}{\hfill \rule{1ex}{1.5ex}\par}
\renewcommand{\[}{$$} \renewcommand{\]}{$$}

\usepackage{tikz}
\usetikzlibrary{positioning,calc,quotes,angles}
\usetikzlibrary{backgrounds}

\usepackage{amsmath}
\usepackage{amssymb,amsfonts}
\usepackage{newtxmath} \usepackage{newtxtext}
\usepackage{mathtools}
\usepackage{fancyvrb,shortvrb}
\usepackage{listings,jlisting,setspace}

\usepackage{wrapfig}
\usepackage{geometry}
\usepackage{subsupscripts}
\DefineShortVerb{\!} 

\numberwithin{equation}{section}
\allowdisplaybreaks
\title{Deformed super brackets on forms of a manifold}

\author{Kentaro Mikami\thanks{Akita University, Japan}\and Tadayoshi
Mizutani\thanks{Saitama University, Japan} }

\AtBeginDocument{\parindent=0pt \setlength{\headheight}{16pt}
  \abovedisplayskip     =0.5\abovedisplayskip
  \abovedisplayshortskip=0.5\abovedisplayshortskip
  \belowdisplayskip     =0.5\belowdisplayskip
  \belowdisplayshortskip=0.5\belowdisplayshortskip}
\renewcommand{\dphi}{d_{\phi}}

\newcommand{\BktT}[2]{\{#1,#2\}^{t}}
\begin{document}
\maketitle
\thispagestyle{plain}

\section{Introduction}
Let \(M\) be a differentiable manifold and \(\tbdl{M}
\) and \(\cbdl{M}\) be the tangent and cotangent bundle of \(M\) (or the
spaces of the sections).   
In addition to the typical example of Lie superalgebra \( \frakg = \sum _{p=1}
^{\dim M} \gaiseki{p} \)  with the Schouten bracket, 
 the space of differential forms 
\( \frakh = \sum_{q=0}^{\dim M} \cgaiseki{q} \) with   
the bracket
\begin{equation}  \BktU{}{\alpha}{\beta} = \parity{a}
d( \alpha \wedge \beta  )\quad \text{where}   \alpha \in \cgaiseki{a}\ 
\text{ and }\
\beta \in \cgaiseki{b} 
\end{equation} 
becomes a Lie superalgebra, where the grading of \( \cgaiseki{a} \) is \( - a -1
\), and is often referred to as \(a'\)    
(cf.\  \cite{Mik:Miz:superForms}).    
The grading of
\( \gaiseki{a} \) is \(  a -1 \), and is also represented by \(a'\).    

There is a notion of deformation of the exterior differentiation  
\(d\) by a 1-form \(\phi\) defined by 
\( d_{t} \alpha = d \alpha +  t \phi \we \alpha \) where \(t\) is a scalar
parameter runs at least \([0,1]\) interval, and   
it is well-known that 
\( d_{t} \circ d_{t} = 0  \) if \(\phi\) is a 1-cocycle.    
It is natural to expect 
\(\dphi\) defines a Lie superalgebra structure, namely  
\begin{equation*}
 \BktT{\alpha}{\beta}=  
\parity{a} d_{t} ( \alpha \we \beta )
= \BktU{}{\alpha}{\beta} +   {\alpha} \we(t\phi)\we {\beta} \end{equation*} 
will be a super bracket for each \(t\).  
Super symmetry holds good. About super Jacobi identity, 
\begin{align*}
& \quad 
\BktT{ \BktT{\alpha}{\beta}}{\gamma} \\ 
& = \Pkt{\Pkt{\alpha}{\beta}}{\gamma} 
+ \parity{a+1} d (  \alpha\wedge\beta \wedge t \phi  \wedge \gamma) 
+ \parity{a} d({\alpha}\wedge{\beta}) \wedge  t \phi \wedge \gamma 
\\ &= \Pkt{\Pkt{\alpha}{\beta}}{\gamma} 
+ \parity{b+c+1}  \alpha\wedge\beta \wedge d \gamma\wedge  t \phi    
+ \parity{b+1}  \alpha\wedge\beta  \wedge  \gamma\wedge d( t \phi)  
\\ &= \Pkt{\Pkt{\alpha}{\beta}}{\gamma} 
+ \parity{b+c+1}  \alpha\wedge\beta \wedge d \gamma\wedge  t \phi    
\quad \text{if  \(\phi\) is closed.}
\\
&\quad  
\mathop{\mathfrak{S}}_{\alpha,\beta,\gamma} 
\parity{a'c'} \BktT{ \BktT{\alpha}{\beta}}{\gamma} 
\\
&=\parity{a'c'}  \Pkt{\Pkt{\alpha}{\beta}}{\gamma} 
+ \parity{a'c'+b+c+1}  \alpha\wedge\beta \wedge d \gamma\wedge  t \phi    
+ \parity{a'c'+ b+1}  \alpha\wedge\beta  \wedge  \gamma\wedge d( t \phi)  
\\&=  
\mathop{\mathfrak{S}}_{\alpha,\beta,\gamma} 
\parity{a c }  d( \alpha) \wedge\beta  \wedge d \gamma \wedge t \phi 
+ 3 
\parity{ac+a+b+c }  \alpha \wedge \beta  \wedge  \gamma \wedge d( t \phi )
\\&=  
\mathop{\mathfrak{S}}_{\alpha,\beta,\gamma} 
\parity{a c }  d( \alpha) \wedge\beta  \wedge d \gamma \wedge t \phi 
\quad \text{if  \(\phi\) is closed.}
\end{align*}

So far, there is no affirmative statement in general setting. 
On the other hand, 
there is a result of deformation of the Schouten bracket by 
	D.~Iglesias  and J.~C.~Marrero 
        in \cite{Igles:Marrero}.  They say 
for a 1-cocycle \(\phi\), 
\begin{equation}
\Sbt{P}{Q}^{\phi} = 
\Sbt{P}{Q} +\parity{p} P(\phi) \wedge (q-1) Q + (p-1)P \wedge Q(\phi) 
\quad 
\text{
where }\quad  P(\phi) = \iota_{\phi} P 
\label{defn:deformedSbt}
\end{equation}
satisfies the axioms of bracket of Lie superalgebra. 
\begin{enumerate}
\item \( \SbtPhi{P}{Q} = - \parity{(p-1)(q-1)} \SbtPhi{Q}{P} \) 
\item \( \SbtPhi{P}{ \SbtPhi{Q}{R} } = 
\SbtPhi{\SbtPhi{P}{Q}}{R} + \parity{(p-1)(q-1)} 
\SbtPhi{Q}{\SbtPhi{P}{R}} \) \; .  
\kmcomment{
\item \( \SbtPhi{P}{ Q \wedge R} = 
\SbtPhi{P}{Q} \wedge {R} + 
\parity{(p-1)q} Q \wedge \SbtPhi{P}{R}
+ \parity{p} P(\phi) \wedge Q \wedge R \)

\item \(\SbtPhi{f}{P} = - P(d_{\phi} f)\) for \(f \in C^{\infty}(M)\)\; . 
}
\end{enumerate}

Inspired by the work above, in this note for a given 1-form \(\phi\), 
we fix properties of a function \(F\) so that 
\begin{equation} \parity{a} d ( \alpha \we \beta ) +  F(a,b)
   {\alpha} \we (t \phi)\we {\beta} 
\quad  \quad ( \alpha \in \cgaiseki{a},  \; 
\beta \in \cgaiseki{b}) \; 
\label{deform:bkt:forms}
   \end{equation}
becomes super bracket for \(t\).   
The function 
\(F\) should be defined on  \( \{ (a,b) \in   \mZ^{2}\mid a+b \leqq \dim M - 1 \}\).


Main results in this note is that there are deformations of two super brackets on the space 
\( \frakh = \sum_{q=0}^{\dim M} \cgaiseki{q} \), and there is a natural
extension to a subalgebra of \(\tbdl{M}\).    

Claim 1:
A deformation of the trivial bracket:
For a given 1-form \(\phi\)  
\begin{equation*} 
\BktU{t,\phi }{ \alpha }{\beta } =  F(a,b)
   {\alpha} \we (t \phi)\we {\beta} 
\quad  \quad ( \alpha \in \cgaiseki{a},  \; 
\beta \in \cgaiseki{b}) \;  
   \end{equation*}
is a super bracket on \(\frakh\) when  \(F\) is a symmetric function.    

Claim 2:
A deformation of the standard bracket:
For a given closed 1-form 
\(\phi\),  
\begin{equation*} \BktU{t,\phi}{\alpha}{\beta} 
=  \parity{a}d({\alpha}\wedge{\beta}) 
+ \tfrac{a+b+2}{2} \alpha \wedge t \phi \wedge  \beta
\quad  \quad ( \alpha \in \cgaiseki{a},  \; 
\beta \in \cgaiseki{b}) \;  
\end{equation*}
is a super bracket on \(\frakh\).

Claim 3: An extension of the 
 deformation of the standard bracket:
For a given closed 1-form \(\phi\), Claim 2 says   
\(\frakh,\BktU{t,\phi}{\cdot}{\cdot } \) are   
Lie superalgebra. Let \(\frakgN{0} ' = \{ X\in\tbdl{M} \mid \Lb{X}\phi = 0 \}
\), which is a subalgebra of \(\tbdl{M}\). 
Then 
\((\frakh,\BktU{t,\phi}{\cdot}{\cdot } ) \oplus\frakgN{0}'\) becomes a
Lie superalgebra naturally by the Lie derivative for each \(t\).

\bigskip

Based on these results, there are many issues to be studied. We would like
to develop homology theory of deformed superalgebra.    
When a Lie group \(G\) acts on \(M\), we get 
\( \sum_{p=1}^{\dim M} \Lambda^{p} \tbdli{M}{G} \) of   
\(G\)-invariant multivector fields, 
and 
\( \sum_{q=0}^{\dim M} \Lambda^{q}\cbdli{M}{G} \) of   
\(G\)-invariant differential forms.  
Since the action preserves the 
Jacobi-Lie bracket, the Schouten bracket is preserved by the action.  
Also the action commutes with the differentiation \(d\), the 
Lie superalgebra  bracket is preserved by the action. 
In short, 
\( \sum_{p=1}^{\dim M} \Lambda^{p} \tbdli{M}{G} \) 
and 
\( \sum_{q=0}^{\dim M} \Lambda^{q}\cbdli{M}{G} \) have Lie superalgebra
structures.    
The simplest case is a Lie group acts on itself.  
\(\ovfrakg =  \sum_{p=1}^{n} \Lambda^{p} \frakg\) where \(\frakg =\)Lie
algebra of \(G\), has a Lie superalgebra structure by the Schouten bracket  
(cf.\cite{Mik:Miz:superLowDim}). 
The differential \(d\) gives  
 a Lie superalgebra structure on 
\(\ovfrakh =  \sum_{p=0}^{n} \Lambda^{p} \frakh\), 
where \( \frakh = \frakg^{*} \) (cf.\  \cite{Mik:Miz:superForms}). 
Concrete and fancy examples are presented from those superalgebras.

\section{Deformation from the trivial bracket}
In this section, we study deformed super bracket of the trivial bracket on
\(\frakh\), namely
\begin{equation} \label{def:def:from:triv}
\BktU{t,\phi }{ \alpha }{\beta } =  F(a,b)
   {\alpha} \we (t \phi)\we {\beta} 
\quad  \quad ( \alpha \in \cgaiseki{a},  \; 
\beta \in \cgaiseki{b}) \; . 
   \end{equation}
Then the symmetric property of \(F\) implies the super symmetric property of
\(
\BktU{t,\phi }{ \cdot  }{\cdot }\)  because of   
\begin{equation*}
\BktU{t,\phi }{ \alpha }{\beta } + 
\parity{ (1+a)(1+b) }
\BktU{t,\phi }{\beta } { \alpha } = (F(a,b) - F(b,a)) \alpha \we (t \phi)
\we \beta \;. 
\end{equation*} 
The super Jacobi identity holds automatically because  
\begin{align*}
\BktU{t,\phi }{ 
\BktU{t,\phi }{ \alpha }{\beta } }
{\gamma} &= 
\BktU{t,\phi }{ 
F(a,b) { \alpha } \we (t \phi) \we {\beta } }
{\gamma} = 
F(a,b) F(a+b+1,c) { \alpha } \we (t \phi) \we {\beta } \we (t\phi) \we
\gamma = 0 \; . 
\end{align*}
\begin{prop}
The bracket \eqref{def:def:from:triv}
is a super bracket on \(\frakh \) when \(F\) is a symmetric function 
on \( \{ (a,b) \in   \mZ^{2}\mid a+b \leqq \dim M - 1 \}\).  
\end{prop}
\begin{remark}
In the proposition above, \(\phi\) is not necessarily closed. 
Like nilpotent subalgebras of a Lie algebra, \( \BktU{t,\phi} 
{ \BktU{t,\phi} {\frakh} {\frakh}} {\frakh} = 0\) holds. 
In contrast,  
\( \parity{a' c'}
\BktU{ }{ 
\BktU{ }{ \alpha  }{\beta  } }{\gamma } = \parity{a' c'} \alpha  \we d\beta \we \gamma 
- \parity{b' a'} \beta  \we d\gamma \we d \alpha  
\) holds  
for \( \BktU{}{\alpha }{\beta } =  \parity{a} d( \alpha  \we \beta  )\). 
\end{remark}

We already know that  the superalgebra  \(\frakh = \sum_{i=0}^{\dim M}
\Lambda^{i} \cbdl{M} \) with the bracket \( \Pkt{\alpha}{\beta} = \parity{a} d
( \alpha\wedge\beta)\) has an extension by  \(\frakgN{0} = \tbdl{M} \) 
through Lie derivative  in \cite{Mik:Miz:superForms}. 
Here we study the deformed superalgebra given by  
the bracket \eqref{def:def:from:triv}
has an extension by a subalgebra of \(\tbdl{M}\)  
through the Lie derivative    
\( \Lb{X} = \iota_{X}\circ d + d \circ \iota_{X}\)  with respect to \(X\). 

Let \(F(a,b)\) be a symmetric function on \(\mZ^{2}\) and
\(\phi\) be a 1-form on \(M\) and 
\begin{equation}
\BktU{t,\phi}{\alpha}{\beta} = F(a,b) \alpha \wedge t \phi \wedge \beta   
\quad  \quad ( \alpha \in \cgaiseki{a},  \; 
\beta \in \cgaiseki{b}) \;. 
\end{equation}
For 1-vector field \(X\), we define \( \BktU{t,\phi} {X}{\alpha} = \Lb{X}\alpha\) and 
 \( \BktU{t,\phi} {\alpha}{X} = - \Lb{X}\alpha\). So the super symmetry holds
 good.  
About super Jacobi identity, we check two cases: one is \(X,Y, \alpha\) and
the other is \(X, \alpha, \beta\).  The first case the super Jacobi identity
is just the formula \( \Lb{\Sbt{X}{Y}} = \Lb{X}\Lb{Y} - \Lb{Y} \Lb{X} \). 
 So the super symmetry holds good.  We treat the other case:  Since 
\begin{align*}   & \quad  
   \BktU{t,\phi}{X}{ \BktU{t,\phi}{\alpha}{\beta}}  
  - \BktU{t,\phi}{ \BktU{t,\phi}{X}{\alpha}}{\beta}  
  - \BktU{t,\phi}{\alpha} { \BktU{t,\phi}{X}{\beta}}  
  \\ &= t F(a,b) \Lb{X} ( \alpha \wedge \phi \wedge \beta )
  - t F(a,b) \Lb{X}\alpha \wedge \phi \wedge \beta
  - t F(a,b) \alpha \wedge \phi \wedge \Lb{X}{\beta}
  \\&= t F(a,b) 
   \alpha \wedge ( \Lb{X} \phi) \wedge \beta
\end{align*}
we see that 
\( \Lb{X} \phi = 0 \) is an efficient condition for super Jacobi identity.    
\begin{prop}
Let \(\frakgN{0}' = \{ X\in \frakgN{0} \mid  
 \Lb{X} \phi = 0 \} \), which is a subalgebra of \(\frakgN{0}\). Then 
 ( \( \frakh \oplus \frakgN{0}', \BktU{t,\phi}{\cdot}{\cdot} \)) is 
 a deformed superalgebra. 

 If \(\phi\) is a 1-cocycle, then 
 ( \( \frakh \oplus \frakgN{0}'', \BktU{t,\phi}{\cdot}{\cdot} \)) is 
 a deformed superalgebra, where  
 \(\frakgN{0}'' = \{ X\in \frakgN{0} \mid  
 \inn{X} \phi = 0 \} \), which is a subalgebra of \(\frakgN{0}'\). 
\end{prop}

Some concrete
example will appear in the tail of the next section.

\section{Deformation from the standard bracket}
It is known in \cite{Mik:Miz:superForms} that 
\( \Pkt{\alpha}{\beta} = \parity{a} d ( \alpha\wedge \beta ) \) defines a
super bracket on \(\sum_{i=0}^{\dim M} \Lambda^{i} \cbdl{M} \). Looking at
the deformed Schouten bracket, the bracket we expect is of form \(
\BktU{t,\phi}{\alpha}{\beta} = \Pkt{\alpha} {\beta} + F(a,b) \alpha \wedge t  \phi \wedge \beta\) for some function \(F\)  
on \( \{ (a,b) \in   \mZ^{2}\mid a+b \leqq \dim M - 1 \}\).  

About super symmetric property, we have  
\begin{equation}
\BktU{t,\phi}{\alpha}{\beta}+ \parity{a' b'}  \BktU{t,\phi}{\beta}{\alpha} = 
( F(a,b) -  F(b,a)) \alpha \wedge t \phi \wedge \beta  
\end{equation}
About Jacobi identity, we see 
\begin{align*} 
&\quad 
\mathop{\mathfrak{S}}_{\alpha,\beta,\gamma} 
\parity{a'c'}\left( \BktU{t,\phi}{ \BktU{t,\phi}{\alpha}{\beta}}{\gamma}  - 
 \Pkt{ \Pkt{\alpha}{\beta}}{\gamma} \right)  
\\
&=  \left( F (1+ b+c,a ) -F ( a,b ) -F ( b,c) -F ( c,a ) +F (1+ c+a,b )  \right) 
 \left( -1 \right) ^{ac+a+c} \left( \alpha \we t \phi\we\beta\we d
 \left( \gamma \right)  \right) 
 \\ & 
 - \left( -1 \right) ^{ac+a+c} \left( F (1+ a+b,c ) +F \left(1+ b+c,a \right) -F ( a,b
 ) -F \left( b,c \right) -F \left( c,a \right)  \right) 
 \left( \alpha \we d \left( \beta \right) \we t \phi \we \gamma \right) 
\\ &
 + \left( F ( a,b ) +F ( b,c ) +F ( c,a ) 
 \right)  \left( -1 \right) ^{ac+a+b+c}
 \left( \alpha \we d \left(t \phi \right) \we \beta \we \gamma \right) 
\\ &
- \left( -1 \right) ^{ac+c}
 \left( F (1+ a+b,c ) -F ( a,b ) -F ( b,c
 ) -F ( c,a ) +F (1+ c+a,b )  \right) 
 \left( d \left( \alpha \right)  \we \beta \we t \phi \we\gamma \right) 
 \;. 
 \end{align*}
Thus, if \(\phi\) is a 1-cocycle, then we get the following  sufficient conditions 
for super Jacobi identity 
\begin{align}
0&= F (1+ b+c,a ) +F (1+ c+a,b ) - F ( a,b ) -F ( b,c ) -F ( c,a )   
\;, \label{eq:x:1}
\\ 0&= F (1+ a+b,c ) +F (1+ b+c,a ) -F ( a,b ) -F ( b,c ) -F ( c,a ) 
\;, \label{eq:x:2}
\\ 0&= F(1+ c+a,b ) + 
F (1+ a+b,c )- F ( a,b ) -F ( b,c ) -F ( c,a ) 
\; . 
\end{align}
The difference \( \eqref{eq:x:1} - \eqref{eq:x:2} =0 \) implies
\( F (1+ c+a,b ) =  F (1+ a+b,c )\).  
Putting \(c=0\), we have 
 \(F (1+ a,b ) =  F (1+ a+b,0 ) \) and \(\ds F(a,0) = F(0,0) ( \frac{a} {2}
 + 1)\), we see that the symmetric function \(F\) 
satisfying the above 3 conditions is
\begin{equation}
F(a,b) = \kappa (a+b+2)\ \ \text{where}\  \kappa \text{ is a constant.} 
\end{equation}

Assume 1-form \(\phi\) is not closed.  Then we get the following sufficient
conditions for super Jacobi identity 
\begin{align}
   0&= F(1+ b+c,a ) + F(1+ c+a,b ) \label{nococycle:a}  
\\ 0&= F(1+ a+b,c ) + F(1+ b+c,a ) \label{nococycle:c}  
\\ 0&= F( a,b ) + F( b,c ) +F( c,a ) \label{nococycle:z}  
\\ 0&= F(1+ c+a,b ) + F(1+ a+b,c ) \label{nococycle:b}  
\\\shortintertext{We conclude a symmetric \(F\) satisfying the 4 above
conditions is trivial as follows. Putting \(c=0\) in \eqref{nococycle:z}}
F(a,b) &=  - G(a) - G(b) \quad \text{where}\ G(a) = F(0,a) = F(a,0)\;.
\\\shortintertext{Applying this expression to 
\eqref{nococycle:z}, we have}
0 &= - 2 ( G(a)+ G(b) + G(c) )\;, \ \text{and so }\; G(a) = 0\;,\ F(a,b) = 0 \;,
\ 
\text{i.e., trivial.} 
\end{align}
We summarize above discussion. 
\begin{thm}
The super symmetry of the bracket 
\eqref{deform:bkt:forms} 
\[\parity{a} d ( \alpha \we \beta ) +  F(a,b)
   {\alpha} \we (t \phi)\we {\beta} \] 
 yields \(F(a,b)\) is a symmetric function, i.e.,  \( F(a,b) = F(b,a) \).  

The super Jacobi identity implies 
if \(\phi\) is not a cocycle, i.e., not exact then \( F(a,b) = 0\).   
If \(\phi\) is a cocycle, i.e., 
if \(\phi\) is exact then the super Jacobi identify yields 
\( F(a,b) = \kappa (a+b+2)\).   
\end{thm}

\kmcomment{
Super Jacobi identity yield 
\begin{align}
0& = 
\mathop{ \mathfrak{S}}_{a,b,c}
\parity{a}F(a,b)- \parity{b}F(b,c+a+1)- \parity{c}F(c,a+b+1)\;.  
\end{align}

\parity{a}F(a,b)+ \parity{b}F(b,c)+ \parity{c}F(c,a)- \parity{b}F(b,c+a+1)- 
\parity{c}F(c,a+b+1)\;.  
kmcomment}

\begin{kmCor}
Let \(\phi\) be an exact 1-form.  
\begin{equation} \BktU{t,\phi}{\alpha}{\beta} 
= \Pkt{\alpha}{\beta} + \tfrac{a+b+2}{2} \alpha \wedge t \phi \wedge  \beta
\label{defn:def:bkt:forms}
\end{equation}
where 
\(\Pkt{\alpha}{\beta} 
= \parity{a}d({\alpha}\wedge{\beta}) \),  
\( \alpha \in \Lambda^{a}\cbdl{M}\) and 
\( \beta \in \Lambda^{b}\cbdl{M}\).  
This  bracket satisfies super symmetry and super Jacobi identity, and the space \(\frakh \)
with this bracket is a Lie superalgebra. 
\end{kmCor}
\begin{remark}
This bracket is a super bracket from Theorem above. 
If one proceed to reconfirm that this bracket is a super bracket,  
it is one page exercise.  Symbol calculus, Maple has a package
\texttt{difforms} and it is helpful to study of properties of this bracket. 
\end{remark}

\kmcomment{
\textbf{Proof:}
For simplicity, we denote \( t \phi\) by \(\phi\) here. 
\begin{align*}
\BktU{t,\phi}{\alpha}{\beta}+ \parity{a' b'}  \BktU{t,\phi}{\beta}{\alpha} &= 
\Pkt{\alpha}{\beta} 
+ \tfrac{a'+b'}{2}\alpha  \wedge  \phi \wedge \beta + \parity{a' b'}  \Pkt{\beta}{\alpha} 
+ \parity{a' b'} \tfrac{a'+b'}{2} \beta \wedge \phi \wedge \alpha \\
&=(\Pkt{\alpha}{\beta}+ \parity{a' b'}  \Pkt{\beta}{\alpha}) + \tfrac{a'+b'}
{2}( \parity{ b+ab+a }+ \parity{a' b'})\beta \wedge \phi \wedge \alpha \\
&= 0  \;.  
\end{align*}
\begin{align*}
\BktU{t,\phi}{ \BktU{t,\phi}{\alpha}{\beta}}{\gamma} 
& = \Pkt{  \BktU{t,\phi}{\alpha}{\beta}} {\gamma} 
+ \tfrac{a'+b'+c'}{2} \BktU{t,\phi}{\alpha}{\beta} \wedge  \phi \wedge  \gamma  
\\ &= 
\Pkt{  \Pkt{\alpha}{\beta} + \tfrac{a'+b'}{2} \alpha \wedge \phi \wedge \beta} {\gamma} 
+ \tfrac{a'+b'+c'}{2}
( \Pkt{\alpha}{\beta} + \tfrac{a'+b'}{2} \alpha \wedge \phi \wedge \beta) \wedge \phi\wedge \gamma 
\\
&= \Pkt{\Pkt{\alpha}{\beta}}{\gamma} + \tfrac{a'+b'}{2} \Pkt{\alpha\wedge\phi\wedge\beta}{\gamma} 
+
\tfrac{a'+b'+c'}{2} 
\parity{a} d({\alpha}\wedge{\beta}) \wedge \phi\wedge \gamma \\
&= \Pkt{\Pkt{\alpha}{\beta}}{\gamma} 
+\tfrac{a'+b'}{2}  \parity{a+1+b} d (  \alpha\wedge\phi \wedge\beta \wedge \gamma) 
+\tfrac{a'+b'+c'}{2}  \parity{a} d({\alpha}\wedge{\beta}) \wedge \phi\wedge \gamma \\
&= \Pkt{\Pkt{\alpha}{\beta}}{\gamma} 
+\tfrac{a'+b'}{2} \parity{1+b'} \phi \wedge d(\alpha\wedge\beta \wedge \gamma) 
+\tfrac{a'+b'+c'}{2} \parity{b'} \phi \wedge d(\alpha\wedge\beta)\wedge \gamma 
\end{align*}
Since 
\(\ds \mathop{\mathfrak{S}}_{\alpha,\beta,\gamma} 
\parity{a'c'} \Pkt{ \Pkt{\alpha}{\beta}}{\gamma} = 0 \), in order to show   
\(\ds 
\mathop{\mathfrak{S}}_{\alpha,\beta,\gamma} 
\parity{a'c'} \BktU{t,\phi}{ \BktU{t,\phi}{\alpha}{\beta}}{\gamma}\) vanishes, we treat
the cycle sum of  
the rest of the last equation after ignoring \( \phi \wedge\). Namely 
\begin{align*} \text{Final target} &= 
\mathop{\mathfrak{S}}_{\alpha,\beta,\gamma} \parity{a' c'} 
\left(\tfrac{a'+b'}{2} \parity{1+b'}  d(\alpha\wedge\beta \wedge \gamma) 
+\tfrac{a'+b'+c'}{2} \parity{b'} d(\alpha\wedge\beta)\wedge \gamma 
\right)
\\
&= \mathop{\mathfrak{S}}_{\alpha,\beta,\gamma} \parity{a' c'} 
\tfrac{a'+b'}{2} \parity{1+b'} \left(d(\alpha\wedge\beta) \wedge \gamma
+ \parity{a'+b'} \alpha\wedge\beta \wedge d \gamma\right) 
\\& 
+ \mathop{\mathfrak{S}}_{\alpha,\beta,\gamma} \parity{a' c'} 
\tfrac{a'+b'+c'}{2} \parity{b'} d(\alpha\wedge\beta)\wedge \gamma 
\\ &= \mathop{\mathfrak{S}}_{\alpha,\beta,\gamma} \parity{a' c'} 
\tfrac{a'+b'}{2} \parity{1+a'}  \alpha\wedge\beta \wedge d \gamma
+ \mathop{\mathfrak{S}}_{\alpha,\beta,\gamma} \parity{a' c'} 
\tfrac{c'}{2} \parity{b'} d(\alpha\wedge\beta)\wedge \gamma 
\\ &= \mathop{\mathfrak{S}}_{\alpha,\beta,\gamma} \parity{a' c'} 
\tfrac{a'+b'}{2} \parity{1+a'}  \alpha\wedge\beta \wedge d \gamma
\\&\quad 
+ \mathop{\mathfrak{S}}_{\alpha,\beta,\gamma} \parity{a' c'} 
\tfrac{c'}{2} \parity{b'}\left((d \alpha) \wedge\beta + \parity{a}\alpha\wedge
 d \beta\right) \wedge \gamma 
\\ &= \mathop{\mathfrak{S}}_{\alpha,\beta,\gamma} \parity{a' c'} 
\tfrac{a'+b'}{2} \parity{1+a'}  \alpha\wedge\beta \wedge d \gamma
\\&\quad 
+ \mathop{\mathfrak{S}}_{\alpha,\beta,\gamma} \parity{a' c'} 
\tfrac{c'}{2} \parity{b'}\left((d \alpha) \wedge\beta  \right) \wedge \gamma 
+ \mathop{\mathfrak{S}}_{\alpha,\beta,\gamma} \parity{a' c'} 
\tfrac{c'}{2} \parity{b'}\left(\parity{a}\alpha\wedge
 d \beta\right) \wedge \gamma 
\\ &= \mathop{\mathfrak{S}}_{\alpha,\beta,\gamma} \parity{a' c'} 
\tfrac{a'+b'}{2} \parity{1+a'}  \alpha\wedge\beta \wedge d \gamma
\\&\quad 
+ \mathop{\mathfrak{S}}_{\alpha,\beta,\gamma} \parity{a' c'} 
\tfrac{c'}{2} \parity{b'+ a'(b'+c')}\beta \wedge \gamma \wedge d \alpha  
+ \mathop{\mathfrak{S}}_{\alpha,\beta,\gamma} \parity{a' c'} 
\tfrac{c'}{2} \parity{b'+ a+b' c+ ac }\gamma\wedge \alpha\wedge d \beta
\\ &= \mathop{\mathfrak{S}}_{\alpha,\beta,\gamma} \parity{a' c'} 
\tfrac{a'+b'}{2} \parity{1+a'}  \alpha\wedge\beta \wedge d \gamma
\\&\quad 
+ \mathop{\mathfrak{S}}_{\alpha,\beta,\gamma} \parity{c' b'} 
\tfrac{b'}{2} \parity{a'+ c'(a'+b')}\alpha  \wedge \beta \wedge d \gamma  
+ \mathop{\mathfrak{S}}_{\alpha,\beta,\gamma} \parity{b' a'} 
\tfrac{a'}{2} \parity{c'+ b+c' a+ ba }\alpha\wedge \beta\wedge d \gamma 
\\ &= \mathop{\mathfrak{S}}_{\alpha,\beta,\gamma} \parity{a' c'} 
\tfrac{a'+b'}{2} \parity{1+a'}  \alpha\wedge\beta \wedge d \gamma
\\& \quad + \mathop{\mathfrak{S}}_{\alpha,\beta,\gamma} \parity{a' c'} 
\tfrac{b'}{2} \parity{a'}\alpha  \wedge \beta \wedge d \gamma  
+ \mathop{\mathfrak{S}}_{\alpha,\beta,\gamma} \parity{a' c'} 
\tfrac{a'}{2} \parity{a'}\alpha\wedge \beta\wedge d \gamma 
\\&= 0 \;. 
\end{align*}
\qed
}

\begin{exam}
Take a 2-dimensional Lie algebra with the Lie bracket relations \(
\Sbt{\yb{1}}{\yb{2}} = \yb{1}\) and the \( \zb{1}, \zb{2}\) is the dual
basis so that  \( d \zb{1} = - \zb{1} \we \zb{2}\) and \( d \zb{2} = 0\). 
The chain complex of weight \(-3\) is given by 
\( \wtedC{1}{-3} = \mR(\zb{1}\we \zb{2})\), 
\( \wtedC{2}{-3} = \mR(\zb{1}\mw 1)+ \mR (\zb{2}\mw 1)\), 
\( \wtedC{3}{-3} = \mR(1\mw 1 \mw 1) = \mR ( \mw^{3} 1 )\). 
We refer to the appendix or \cite{Mik:Miz:superForms} about odd notations. 
Fix \(\phi = \zb{2}\).  \(\pdel (\zb{1}\we \zb{2}) = 0\) for 1-chain.   
\begin{align*}
 \pdel (\zb{j}\mw 1 )& = \BktU{t, \zb{2}}{\zb{j}}{1}
= \begin{cases*}
c_{1} \delta_{j}^{1} t \zb{1} \we \zb{2} & \(c_{1}\) is constant,  deform of trivial \\
\delta_{j}^{1} ( 1+ \frac{3}{2}t )
\zb{1} \we \zb{2} & deform of standard \end{cases*} 
\\
\pdel ( \mw ^{3} 1 ) &= \tbinom{3}{2}  \BktU{t, \zb{2}}{1}{1} \mw 1
= \begin{cases*}
c_{0} t \zb{2}\mw 1 & \(c_{0}\) is constant,  deform of trivial \\
 t \zb{2} \mw 1 & deform of standard \end{cases*} 
 \end{align*}   
We summarize the kernel dimensions and Betti numbers as follows, we assume 
\( c_{0} c_{1} \ne 0\)).

\renewcommand{\arraystretch}{0.8}
\[
\begin{array} {c | *{3}{c}}
\text{trivial} & 1 & 2 & 3 \\\hline
\dim & 1 & 2 & 1 \\\hline
\ker \dim & 1 & 1 + \delta_{t}^0 &  \delta_{t}^0 \\ 
\text{Betti} &  \delta_{t}^0 & 2 \delta_{t}^0 &  \delta_{t}^0 \end{array} 
\hspace{20mm}
\begin{array} {c | *{3}{c}}
\text{standard} & 1 & 2 & 3 \\\hline
\dim & 1 & 2 & 1 \\\hline
\ker \dim & 1 & 1 + \delta_{2+3t}^0 &  \delta_{t}^0 \\ 
\text{Betti} &  \delta_{2+3t}^0 & \delta_{2+3t}^0 +  \delta_{t}^0 &  \delta_{t}^0 
\end{array} 
\]
\renewcommand{\arraystretch}{1.0}

\label{exam:w3:dim2}
\end{exam}

\subsection{An extension}

We mentioned before that the superalgebra  \(\frakh = \sum_{i=0}^{\dim M}
\Lambda^{i} \cbdl{M} \) with the bracket \( \Pkt{\alpha}{\beta} = \parity{a}
d ( \alpha\wedge\beta)\) has an extension by  \(\tbdl{M} \) through Lie
derivation in \cite{Mik:Miz:superForms}.  Here we study two  deformed
superalgebras  have an extension by a subalgebra of  \( \frakgN{0} = \tbdl{M}
\) through the Lie derivative    \( \Lb{X} = \iota_{X}\circ d + d \circ
\iota_{X}\)  with respect to \(X\).

The superalgebra  \(\frakh \) has the deformed bracket \( \BktU{t, \phi}
{\alpha}{\beta} = \parity{a} d ( \alpha\wedge\beta)+ \frac{a+b+2}{2} \alpha
\wedge t \phi \wedge \beta \), where \(\phi\) is a 1-cocycle.    

Let \( \Ekt{\cdot}{\cdot} \) be a
candidate of superbracket on \(\frakh \oplus \tbdl{M}\), i.e.,   
\( \Ekt{\alpha}{\beta} = \BktU{t,\phi}{\alpha}{\beta} \) ,  
\( \Ekt{X}{Y} = \Sbt{X}{Y} \),  
\( \Ekt{X}{\beta} = - \Ekt{\beta}{X} = \Lb{X}{\beta} \) for forms \(\alpha,
\beta\) and 1-vectors \(X,Y\).  
We have to check super Jacobi identity for two cases. Again we abbreviate
\(t\phi\) by \(\phi\).  One case is all right as
following. 
\begin{align*}
&\quad \Ekt{ \Ekt{X}{Y} }{\alpha}
+\Ekt{ \Ekt{Y}{\alpha} }{X}
+\Ekt{ \Ekt{\alpha}{X} }{Y}
\\=& \Ekt{ \Sbt{X}{Y} }{ \alpha } - \Lb{X} \Lb{Y} \alpha  + \Lb{Y} \Lb{X} \alpha 
=(\Lb{ \Sbt{X}{Y} }  - \Lb{X} \Lb{Y} + \Lb{Y} \Lb{X}) \alpha 
= 0 \; . 
\end{align*}
We try the other.   
\begin{align*}
&\quad \Ekt{X}{\Ekt{\alpha}{\beta} }
+\Ekt{\alpha}{ \Ekt{\beta}{X } }
+\parity{b'a'}\Ekt{\beta}{\Ekt{X}{\alpha} }
\\ =& \Lb{X} \Ekt{\alpha}{\beta} 
 +  \parity{a} d( \alpha \wedge \Ekt{\beta}{X} )
 + \frac{a'+b'}{2} \alpha \wedge \phi \wedge \Ekt{ \beta}{X}
\\
& + \parity{1+a' b'} \left(
\parity{b'} d ( \beta \wedge  \Lb{X}{\alpha}  )
- \frac{a'+b'}{2} \beta\wedge \phi \wedge  \Lb{X}{\alpha}  
\right)
\\ =& \Lb{X}( \parity{a} d( \alpha\wedge\beta ) + \frac{a'+b'}{2} \alpha\wedge
\phi\wedge \beta)  
 +  \parity{a+1} d( \alpha \wedge \Lb{X}\beta  ) 
 \\& 
- \frac{a'+b'}{2} \alpha \wedge
\phi \wedge \Lb{X}\beta 
 + \parity{1+a' b'} \left(
\parity{b'} d ( \beta \wedge  \Lb{X}{\alpha}   )
- \frac{a'+b'}{2} \beta\wedge \phi \wedge  \Lb{X}{\alpha}   
\right)
\\ =&  \parity{a} d \Lb{X} ( \alpha\wedge\beta )    
+ \frac{a'+b'}{2}\Lb{X}(\alpha\wedge \phi\wedge \beta)  
\\ & +  \parity{a+1} d( \alpha \wedge\Lb{X}\beta  ) 
- \frac{a'+b'}{2} \alpha \wedge
\phi \wedge (\Lb{X}\beta ) 
\\
& + \parity{1+a  b'} 
 d ( \beta \wedge  \Lb{X}{\alpha}   ) 
+ \parity{a' b'}\frac{a'+b'}{2} \beta\wedge \phi \wedge 
( \Lb{X}{\alpha} + \inner{X}{00} \alpha ) 
\\ =& 
+ \frac{a'+b'}{2}\Lb{X}(\alpha\wedge \phi\wedge \beta)  
 + \parity{a' b'}\frac{a'+b'}{2} \beta\wedge \phi \wedge 
 \Lb{X}{\alpha}   
- \frac{a'+b'}{2} \alpha \wedge
\phi \wedge (\Lb{X}\beta ) 
\\ &
+ \parity{a} d \Lb{X}( \alpha\wedge\beta  )  
+  \parity{a+1} d( \alpha \wedge\Lb{X}\beta ) 
 + \parity{a+1 } 
 d ( \Lb{X}{\alpha} \wedge \beta   ) 
\\ =& 
+ \frac{a'+b'}{2}\Lb{X}(\alpha\wedge \phi\wedge \beta)  
+ \parity{a' b'}\frac{a'+b'}{2}(\beta\wedge \phi \wedge 
 \Lb{X}{\alpha} ) 
- \frac{a'+b'}{2}(\alpha \wedge \phi \wedge \Lb{X}\beta ) 
\\ =& 
+ \frac{a'+b'}{2}\left(\Lb{X}(\alpha\wedge \phi\wedge \beta)  
-  \Lb{X}{\alpha}\wedge \phi \wedge \beta  
- \alpha \wedge \phi \wedge \Lb{X}\beta \right)  
\\ =& 
+ \frac{a'+b'}{2}\left(\alpha\wedge \Lb{X}(\phi) \wedge \beta  
\right)  
\end{align*}
This vanishes if \(\Lb{X}\phi=0\). We see that 
\(\{ X\in \tbdl{M} \mid  \Lb{X}\phi=0\}\) forms  a \(\mR\) subalgebra of \(\tbdl{M}\). 
Thus we have the following result. 

\begin{thm}

The superalgebra  \((\frakh,\BktU{t,\phi}{\cdot}{\cdot}) \) allows an
extension by \(\frakgN{0}' = \{ X\in \tbdl{M} \mid  \Lb{X}\phi=0\}\),
namely, \( \frakh \oplus \frakgN{0}' \) becomes a Lie superalgebra extension
of \(\frakh\).  
\end{thm}

\begin{exam}
We extend \(\frakh\) in Example \ref{exam:w3:dim2} by \( \frakgN{0}\) 
and show the \((- 3)\)-weighted chain complex. 
We denote \( \yb{1}\mw \yb{2}\) by \(U\), and   
 \( \zb{1}\we \zb{2}\) by \(V \).  
\[ 
 \begin{array}[t] {c|*{5}c | c}
\wtedCR{m}{3}  & 1 & 2 & 3 & 4 & 5 & \text{ if }\\\hline
\dim & 1 & 4 & 6 & 4 & 1
\\\hline
\text{basis} & V & 
{ \renewcommand{\arraystretch}{0.7}
\begin{array}{c}
\zb{j}\mw 1 \\ \yb{i}\mw V
\end{array} } & 
{ \renewcommand{\arraystretch}{0.7}
\begin{array} {c}
1^{3} \\
\yb{i}\mw \zb{j} \mw 1 \\ U \mw V
\end{array} } & 
{ \renewcommand{\arraystretch}{0.7}
\begin{array}{c}
\yb{i}\mw 1^{3} \\ U \mw \zb{j} \mw 1  
\end{array} } & U \mw 1^{3} 
\\\hline
\ker\dim & 1 & 3 & \begin{array}{c} 4 \\ 3 \end{array}  
                 & \begin{array}{c} 2 \\ 1 \end{array} 
 & 0 & \begin{array}{c} t(1+3t/2) = 0 \\ t(1+3t/2) \ne  0 \end{array} \\
 \hline
\text{Betti} & 0 & \begin{array}{c} 1 \\ 0 \end{array}  
                 & \begin{array}{c} 2 \\ 0 \end{array} 
                 & \begin{array}{c} 1 \\ 0 \end{array} 
 & 0 & \begin{array}{c} t(1+3t/2) = 0 \\ t(1+3t/2) \ne  0 \end{array} \\
\end{array} \]

By the direct computation, we see that 
the boundary image is spanned as follows. 
\begin{align*}
\pdel \wtedC{1}{-3} & = \{ 0 \}\;, \qquad  
\pdel \wtedC{2}{-3} = \{ (1+3t/2) V, V \}\\ 
\pdel \wtedC{3}{-3} &= \{ 3t \zb{2}\mw 1, 
 \zb{2}\mw 1 + (1+3t/2) \yb{1} \mw V, 
 \zb{1}\mw 1 - (1+3t/2) \yb{2} \mw V \}\;, \\ 
\pdel \wtedC{4}{-3} & = \{ t \yb{1} \mw \zb{2}\mw 1,  
t \yb{2} \mw \zb{2}\mw 1, - \yb{2} \mw \zb{2}\mw 1 + (1+3t/2) \yb{1} \mw \yb{2} \mw V, 
  \yb{1} \mw \zb{2}\mw 1\}\;, \\ 
\pdel \wtedC{5}{-3} & = \{  \yb{1} \mw^{3}  1 + 3t  \yb{1} \mw \yb{2} \mw
\zb{2} \mw 1 \}\; .
\end{align*}   
Thus, the kernel dimensions of \(\pdel\) for m=1,3,5 are 1,3,0 and those of m = 3,4 are
4,2 if \(t(1+3t/2) = 0\) else  
3,1. Finally the Betti numbers are 
0,1,2,1, 0 if \(t(1+3t/2) = 0\) else  0,0,0,0,0 .  

To get the weighted homology groups of Lie algebras, even of 2-dimensional, 
we need hard work. We prepare reporting the general weighted homology
groups of 2-dimensional Lie algebra. 
\end{exam}

In \cite{Mik:Miz:super2} and \cite{Mik:Miz:super3}, we introduced double
weight for the algebra of homogeneous polynomial coefficient multi-vector
fields on \(\mR^{n}\).  By the similar way, we get examples of double
weighted super algebras of homogeneous polynomial coefficient forms and
1-vector fields on  \(\mR^{n}\) in \cite{Mik:Miz:superForms}. 
It may be interesting to 
study those deformed  double weighted superalgebras and their homology
groups.  


\appendix
\def\thesection{Appendix: \Alph{section}}

\section{Quick review of the homology groups of Lie superalgebra}

Let \( \frakg = \sum_{i\in\mZ} \frakgN{i} \) be  a Lie superalgebra.
From super symmetry 
\( \Sbt{X}{Y} + \parity{ x y } \Sbt{Y}{X} = 0 \) for \(  
 X \in \frakgN{x}, Y \in \frakgN{y}  \),  
\(m\)-th chain space is given by  \(\myCS{m} = \otimes ^{m} \frakg /\text{Ideal of} 
( {X} \otimes {Y} + \parity{ x y } {Y} \otimes {X} )\).   
We denote the class of \( A_{1}  \otimes   \cdots  \otimes   A_{p}\) by   
\( A_{1}  \mw  \cdots  \mw  A_{p}\). 
Let \( \widetilde{A} = A_{1}  \mw  \cdots  \mw  A_{p}\) and  
 \( \widetilde{B} = B_{1}  \mw  \cdots  \mw  B_{q}\).    
The boundary operator \(\ds \pdel_{} :\myCS{m}\to \myCS{m-1}\)
called (\textit{boundary homomorphism})  is defined by 
\begin{align} 
 \pdel ( Y_{1}\mywedge \cdots \mywedge Y_{m}  ) 
&= 
\sum_{i<j} (-1)^{ i-1 + y_{i}(\mathop{\sum}_{i< s<j} y_{s}) } 
Y_{1} \mywedge \cdots \widehat{ Y_{i} } \cdots \mywedge 
\underbrace{\Sbt{Y_{i}}{Y_{j}}}_{j} \mywedge \cdots  \mywedge Y_{m}
\label{triv:1}
\end{align}
for a decomposable element, where \(\ds y_{i}\) is the degree of homogeneous element
 \(Y_{i}\), i.e., \(\ds Y_{i} \in \frakg_{y_{i}} \). 
It is clear that \(\ds
\pdel\circ \pdel = 0\) and we have the homology groups 
\(\ds 
\myHom{m}(\frakg, \mR) =  \ker(\pdel : \myCS{m} \rightarrow
\myCS{m-1})/ \pdel ( \myCS{m+1} )
\).  

We say 
 a non-zero \(m\)-th decomposable element 
 \( Y_{1}\mywedge \cdots \mywedge Y_{m}  \) has the weight  
 \( \sum_{i=1}^{m} \yb{i} \) where  
\(\ds Y_{i} \in \frakgN{\yb{i}} \). 
The weight is preserved by \(\pdel\), i.e., 
\( \pdel ( \wtedC{m}{w} ) \subset \wtedC{m-1}{w}\) where  
 \(\wtedC{m}{w} =  \) the subspace of \(w\)-weighted  \(m\)-th chains and 
 we have the weighted homology groups.   
\begin{align} 
\shortintertext{
If all $y_{i}$ are even in \eqref{triv:1}, then }
 \pdel ( Y_{1}\mywedge \cdots \mywedge Y_{m}  ) 
& = -  \sum_{i<j} (-1)^{ i+j }
 \Sbt{Y_{i}}{Y_{j}} \mywedge 
 Y_{1} \mywedge \cdots\widehat{Y_{i}}\cdots
\widehat{Y_{j}}\cdots \mywedge Y_{m} 
\;. \label{all:even}
\\
\shortintertext{
If all $y_{i}$ are odd in \eqref{triv:1}, then }
 \pdel ( Y_{1}\mywedge \cdots \mywedge Y_{m}  ) 
& =  \sum_{i<j} 
 \Sbt{Y_{i}}{Y_{j}} \mywedge 
 Y_{1} \mywedge \cdots\widehat{Y_{i}}\cdots
\widehat{Y_{j}}\cdots \mywedge Y_{m}   \;.  \label{all:odd}
\end{align}
\begin{defn}
Let 
\(A =A_{1} \mywedge \cdots \mywedge A_{\bar{a}} \) 
(\( A_{i}\in \frakg_{a_{i}}\))    
and \(B =B_{1} \mywedge \cdots \mywedge B_{\bar{b}} \) 
( \( B_{j}\in \frakg_{b_{j}}\)).    
Define 
\begin{equation}
\SbtES{A}{B} = 
\pdel ( A\mywedge B ) - (\pdel A)\mywedge B - (-1)^{\bar{a}} A
\mywedge \pdel B \;.
\label{bdary:bunkai}
\end{equation}
\end{defn}
It satisfies 
\begin{align}
\SbtES{A}{B} &= \sum_{i,j} (-1)^{i+ a_{i} \sum\limits_{s>i}
a_{s}+ j+ b_{j}(1+ \sum\limits_{s=1}^{j} b_{s})}
A_{1} \mywedge 
\cdots \widehat{ A_{i} }\cdots A_{\bar{a}} \mywedge \Sbt{A_{i}}{B_{j}} 
\mywedge
B_{1} \mywedge 
\widehat{ B_{j} }
\cdots \mywedge B_{\bar{b}}\;.  \label{defn:super:Schouten}
\\
\noalign{ If all \(a_{i}\) are even and all \(b_{j}\) are odd
in \eqref{defn:super:Schouten}, then} 
\SbtES{A}{B}  
& =   \sum_{i,j} (-1)^{i+1} A_{1}\mywedge\cdots\widehat{ A_{i} }\cdots 
A_{\bar{a}} 
\mywedge \Sbt{A_{i}}{B_{j}} \mywedge
B_{1} \mywedge 
\cdots \widehat{ B_{j} }\cdots \mywedge B_{\bar{b}}\;. 
\\\shortintertext{
Let 
\(C =C_{1} \mywedge \cdots \mywedge C_{\bar{c}} \) 
(\( C_{k}\in \frakg_{c_{k}}\)) with \( c_{k}\) are all even. Then     
}
\SbtES{A}{C\mw B}  
& =   \sum_{i,k} (-1)^{i+k} A_{1}\mywedge\cdots\widehat{ A_{i} }\cdots A_{\bar{a}} 
\mywedge \Sbt{A_{i}}{C_{k}} \mywedge C_{1} \mywedge 
\cdots \widehat{ C_{k} }\cdots \mywedge C_{\bar{c}} \mw B 
\\&\quad 
+ \sum_{i,j} (-1)^{i+1} A_{1}\mywedge\cdots\widehat{ A_{i} }\cdots A_{\bar{a}} 
\mw C 
\mywedge \Sbt{A_{i}}{B_{j}} \mywedge B_{1} \mywedge 
\cdots \widehat{ B_{j} }\cdots \mywedge B_{\bar{b}} \;.
\notag
\\
& =  \sum_{i,k} (-1)^{i+1} A_{1}\mywedge\cdots\widehat{ A_{i} }\cdots A_{\bar{a}} 
\mywedge C_{1} \mywedge \cdots \mw \Sbt{A_{i}}{C_{k}} \mw 
\cdots \mywedge C_{\bar{c}} \mw B 
\\&\quad 
+ \sum_{i,j} (-1)^{i+1} A_{1}\mywedge\cdots\widehat{ A_{i} }\cdots A_{\bar{a}} 
\mw C 
\mywedge \Sbt{A_{i}}{B_{j}} \mywedge B_{1} \mywedge 
\cdots \widehat{ B_{j} }\cdots \mywedge B_{\bar{b}} \;.
\notag
\end{align}


\nocite{Mik:Miz:super2}
\nocite{Mik:Miz:super3}

\bibliographystyle{plain}
\bibliography{km_refs}

\end{document}